\documentclass[11pt,reqno,oneside]{amsart}
 \usepackage{latexsym}
 \usepackage{amssymb}
 \usepackage{amsfonts}
\usepackage{graphics}
\usepackage{mathcomp}
 \author[Poncin]{N. Poncin*}
 \author[Radoux]{F. Radoux**}
 \author[Wolak]{R. Wolak***}

\usepackage{amsmath}\usepackage{epsf,amsfonts,amsthm}\usepackage{amscd,mathptmx,amssymb}
\usepackage{xcolor,epic,eepic}\usepackage{epsfig}
\usepackage{fontenc,indentfirst, delarray,amsfonts,amsmath,amssymb}
\usepackage{rotating}
\usepackage[T1]{fontenc}

\thanks{* University of Luxembourg, Campus Limpertsberg, Mathematics Research Unit,
162A, avenue de la Fa\"iencerie, L-1511 Luxembourg
City, Grand-Duchy of Luxembourg, E-mail: norbert.poncin@uni.lu.\\
** University of Li\`ege, Institute of Mathematics, Grande Traverse, 12 - B37, B-4000
 Li\`ege, Belgium, E-mail: Fabian.Radoux@ulg.ac.be.\\
*** Jagiellonian
University, ulica Reymonta 4 30-059 Krakow, Poland, E-mail:
Robert.Wolak@im. uj.edu.pl.}
\date{\today}
\title[Singular Quantization]{Equivariant quantization of orbifolds}

\newtheorem{lem}{Lemma}
\newtheorem{thm}[lem]{Theorem}
\newtheorem{prop}[lem]{Proposition}

\theoremstyle{remark}

\theoremstyle{definition}
\newtheorem{defi}{Definition}

\newcommand{\R}{\mathbb{R}}

\newcommand{\N}{\mathbb{N}}

\newcommand{\Ci}{C^{\infty}}

\newcommand{\nc}{\newcommand}
\nc{\rnc}{\renewcommand} 
\linethickness{.15ex}

\begin{document}
\begin{abstract}
Equivariant quantization is a new theory that highlights the role
of symmetries in the relationship between classical and quantum
dynamical systems. These symmetries are also one of the reasons
for the recent interest in quantization of singular spaces,
orbifolds, stratified spaces... In this work, we prove existence
of an equivariant quantization for orbifolds. Our construction
combines an appropriate desingularization of any Riemannian
orbifold by a foliated smooth manifold, with the foliated
equivariant quantization that we built in [Poncin N, Radoux F, Wolak R, {\it A first approximation for
quantization of singular spaces},
J. Geom. Phys., {\bf 59} (4) (2009), pp 503-518]. Further,
we suggest definitions of the common geometric objects on
orbifolds, which capture the nature of these spaces and guarantee,
together with the properties of the mentioned foliated resolution,
the needed correspondences between singular objects of the
orbifold and the respective foliated objects of its
desingularization.
\end{abstract}
\maketitle
{\bf{Subject classification}} : Real and complex differential geometry. 
\section{Introduction}

{\sl Equivariant quantization}, see \cite{LMT}, \cite{LO},
\cite{DLO}, \cite{PL}, \cite{BM}, \cite{DO}, \cite{BHMP},
\cite{BM2}... is the fruit of a recent research program that aimed
at a complete and unambiguous geometric characterization of
quantization. The procedure highlights the primary role of
symmetries in the relationship between classical and quantum
dynamical systems. One of the major achievements of equivariant
quantization is the understanding that a fixed $G$-structure of
the configuration space of a mechanical system guarantees
existence and uniqueness of a $G$-equivariant quantization.
Roughly and more generally, an equivariant, or better, a natural
quantization of a smooth manifold $M$ is a vector space
isomorphism
$$Q[\nabla]:\mathrm{Pol}(T^*M)\ni s\to
Q[\nabla](s)\in\mathcal{D}(M)$$ that maps a smooth function
$s\in\mathrm{Pol}(T^*M)$ of the phase space $T^*M$, which is
polynomial along the fibers, to a differential operator
$Q[\nabla](s)\in\mathcal {D}(M)$ that acts on functions $f\in
C^{\infty}(M)$ of the configuration space $M$. The quantization
map $Q[\nabla]$ depends on the projective class $[\nabla]$ of an
arbitrary torsionless connection $\nabla$ of $M$, and it is
equivariant with respect to the action of local diffeomorphisms
$\phi$ of $M$, i.e.
$$Q[\phi^*\nabla](\phi^*s)(\phi^*f)=\phi^*(Q[\nabla](s)(f)),$$
$\forall s\in\mathrm{Pol}(T^*M),\forall f\in C^{\infty}(M).$ Such
natural and projectively invariant quantizations, or simply
equivariant quantizations, were investigated in several works, see
e.g. \cite{MB}, \cite{MR}, \cite{SH}.\medskip

On the other hand, {\sl quantization of singular spaces}, see e.g.
\cite{BHPSingReduDefoQuan}, \cite{JHQuantReduc},
\cite{JHSingQuant}, \cite{JRSQuantStrat},
\cite{HLSingQuanCommRedu}, \cite{PflDefoQuanSympOrbi}... is an
upcoming topic in Mathematical Physics, in particular in view of
the interest of reduction for complex systems with symmetries.
More precisely, if a symmetry group acts on the phase space or the
configuration space of a general system, the quotient space is
usually a singular space, an orbifold or a stratified space... The
challenge consists in the quest for a quantization procedure of
such singular spaces that in addition commutes with
reduction.\medskip

It is now quite natural to ask which aspects of the new theory of
equivariant quantization -- that was recently extended from vector
spaces to smooth manifolds -- hold true for certain singular
spaces. The main result of this work is the proof of existence of
{\sl equivariant quantization for orbifolds}.\medskip

A first difficulty of the attempt to construct an equivariant
quantization on a singular space, is the proper definition of the
actors in equivariant quantization -- functions, differential
operators, symbols, vector fields, differential forms,
connections... -- for this space. Even in the case of orbifolds no
universally accepted definitions can be found in literature.
Morevoer, geometric and algebraic definitions do not always
coincide as in the classical context. Our method is based on the
resolution of orbifolds proposed in \cite{Girbau}. More precisely,
we combine this desingularization technique, which allows
identifying any Riemannian orbifold $V$ with the leaf space of a
foliated smooth manifold $(\tilde{V},\mathcal{F})$, with the
foliated equivariant quantization that we constructed in
\cite{PoRaWo}, to build a singular equivariant quantization of
orbifolds. To realize this idea, meaningful definitions, which not
only capture the nature of orbifolds but ensure simultaneously
that singular objects of $V$ are in 1-to-1 correspondence with the
respective foliated objects of $(\tilde{V},\mathcal F)$, are
needed. We show that the chosen foliated resolution of orbifolds
has exactly the properties that are necessary for this kind of
relationship.\medskip

The paper is organized as follows. In the second section, we
recall the definitions of foliated objects and of a foliated
equivariant quantization. In the third, we detail our geometric
definitions of singular objects on orbifolds and study their
relevant properties for the singular equivariant quantization
problem. We describe and further investigate, in Section 4, the
foliated desingularization of a Riemannian orbifold, putting
special emphasis on aspects that are of importance for the
mentioned appropriate correspondence between foliated and singular
objects. The last section deals with existence and the explicit
construction of a singular equivariant quantization of Riemannian
orbifolds.

\section{Foliated quantization}

In the sequel, $(M,\mathcal{F})$ denotes an $n$-dimensional smooth manifold endowed with a regular foliation $\mathcal{F}$ of dimension $p$ and codimension $q=n-p$. Moreover, $U$ is an open set of $(M,\mathcal{F})$.

Let us first recall the definitions of the foliated objects and of
the foliated natural and projectively invariant quantization given
in \cite{PoRaWo} :

\begin{defi}
A {\sl foliated function} $f$ on $U$ is a smooth function $f\in C^{\infty}(U)$ such that $f$ is constant along the connected components of the traces of the leaves in $U$. In other words, if $(V,(x,y))$ is a system of adapted coordinates such that $V\cap U\neq\emptyset$, the local form of $f$ on $U\cap V$ depends only on the transverse coordinates $y$. \end{defi}

We denote by $C^{\infty}(U,\mathcal{F})$ the algebra of all
foliated functions of $(U,\mathcal{F})$.

\begin{defi}
A {\sl foliated differential operator} $D$ of order $k\in\N$ of $U$ is an endomorphism of the space $\Ci(U,\mathcal{F})$ of foliated functions, which reads in any system $(V,(x^1,\ldots,x^p,y^1,\ldots,y^q))$ of adapted coordinates in the following way:
$$D|_{U\cap V}= \sum_{\vert \alpha \vert \le k}
D_{\alpha}\,\partial_{y^{1}}^{\alpha^1}\ldots\partial_{y^q}^{\alpha^q},$$
where the coefficients $D_{\alpha}\in C^{\infty}(U\cap V,\mathcal{F})$
are locally defined foliated functions and where $k$ is
independent of the considered chart.
\end{defi}

We denote by $\mathcal{D}^k(U,\mathcal{F})$ the
$\Ci(U,\mathcal{F})$-module of all $k$-th order foliated
differential operators of $(U,\mathcal{F})$ and set
$$\mathcal{D}(U,\mathcal{F}):= \cup_{k\in\mathbb{N}}\mathcal{D}^k(U,\mathcal{F}).$$

The graded space $\mathcal{S}(U,\mathcal{F})$ associated with the filtered space
$\mathcal{D}(U,\mathcal{F})$,
$$\mathcal{S}(U,\mathcal{F}):=\oplus_{k\in\mathbb{N}}\mathcal{S}^k(U,\mathcal{F}):=\oplus_{k\in\mathbb{N}}\mathcal{D}^k(U,\mathcal{F})/\mathcal{D}^{k-1}
(U,\mathcal{F}),$$ is the space of {\sl foliated symbols}. The
$k$-th order symbol of a $k$-th order foliated differential
operator $D$ is then simply its class $\sigma_k(D)$ in the $k$-th
term of the symbol space. The principal symbol $[D]$ of $D$ is the
symbol $\sigma_k(D)$ with the lowest possible $k$.

\begin{defi}
An {\sl adapted vector field} of $U$ is a vector
field $X\in\mathrm{Vect}(\sl{U})$ such that
$[X,Y]\in\Gamma(T\mathcal {F})$, for all $Y\in\Gamma(T\mathcal
{F})$.
\end{defi}
The space $\mathrm{Vect}_{\mathcal{F}}(\sl{U})$ of adapted
vector fields is obviously a Lie subalgebra of the Lie algebra
$\mathrm{Vect}(\sl{U})$ and the space $\Gamma(T\mathcal{F})$ of tangent vector
fields is an ideal of $\mathrm{Vect}_{\mathcal{F}}(U)$.
\begin{defi}
The quotient algebra
$\mathrm{Vect}(\sl{U},\mathcal{F}):=\mathrm{Vect}_{\mathcal{F}}(\textit{U})/\mathrm{\Gamma}(\textsl{T}\mathcal{F})$
is the Lie algebra of {\sl foliated vector fields}.
\end{defi}

The space $\mathrm{Vect}(U,\mathcal{F})$ is also a
$\Ci(U,\mathcal{F})$-module that acts naturally on $\Ci(U,$
$\mathcal{F})$.

\begin{prop}
The space $\mathrm{Vect}(\sl{U},\mathcal{F})$ is isomorphic to the
space $\mathcal{S}^1(U,\mathcal{F})$.
\end{prop}
\begin{proof}
See \cite{PoRaWo}.
\end{proof}
\begin{defi}
A {\sl foliated differential $1$-form} of $U$ is a
differential $1$-form $\theta$ of $U$ such that
$\mathrm{i}_{\sl{Y}}\theta=\mathrm{i}_{\sl{Y}}\mathrm{d}\theta=0,$
for all $Y\in\Gamma(T\mathcal{F})$.
\end{defi}

We denote by $\Omega^1(U,\mathcal{F})$ the space of all foliated
differential $1$-forms of $U$. The interior product
of a foliated $1$-form with a foliated vector field is a foliated
function.

\begin{defi} A {\sl foliated torsion-free connection} of $U$ is a
bilinear map
$\nabla(\mathcal{F}):\mathrm{Vect}(U,\mathcal{F})\times\mathrm{Vect}(U,\mathcal{F})\to\mathrm{Vect}(U,\mathcal{F})$
such
that, for all $f\in C^{\infty}(U,\mathcal{F})$ and all
$[X],[Y]\in\mathrm{Vect}(U,\mathcal{F})$, the following conditions
are satisfied:
\begin{itemize}\item
$\nabla(\mathcal{F})_{f[X]}[Y]=f\nabla(\mathcal{F})_{[X]}[Y]$,\item
$\nabla(\mathcal{F})_{[X]}(f[Y])=([X].f)[Y]+f\nabla(\mathcal{F})_{[X]}[Y],$\item
$\nabla(\mathcal{F})_{[X]}[Y]=\nabla(\mathcal{F})_{[Y]}[X]+[[X],[Y]].$
\end{itemize}
\end{defi}

We denote by $\mathcal{C}(U,\mathcal{F})$ the affine space of
torsion-free foliated connections of $U$.

\begin{defi}
Two foliated connections $\nabla(\mathcal{F})$ and
$\nabla'(\mathcal{F})$ of $U$
are {\sl projectively equivalent} if and only if there is a
foliated $1$-form $\theta\in\Omega^1(U,\mathcal{F})$ such that,
for all $[X],[Y]\in\mathrm{Vect}(U,\mathcal{F})$, one has
$$\nabla'(\mathcal{F})_{[X]}[Y]-\nabla(\mathcal{F})_{[X]}[Y]=\theta([X])[Y]+\theta([Y])[X].$$
\end{defi}

\begin{defi} A {\sl foliated local diffeomorphism} between two
foliated manifolds $(M,\mathcal{F})$ and $(M',\mathcal{F'})$ is a
smooth mapping $\Phi:M\to M'$ that is locally a diffeomorphism and
maps any leaf of $\mathcal{F}$ into a leaf of $\mathcal F'$.
\end{defi}

\begin{defi}
A {\sl foliated natural and projectively invariant quantization}
is a map
$${\mathcal
Q}(\mathcal{F}):\mathcal{C}(M,\mathcal{F})\times\mathcal{S}(M,\mathcal{F})\to\mathcal{D}(M,\mathcal{F}),$$
which is defined for any foliated manifold $(M,\mathcal{F})$ and
has the following properties:
\begin{itemize}
\item ${\mathcal Q}(\mathcal{F})(\nabla(\mathcal{F}))$ is a linear
bijection between $\mathcal{S}(M,\mathcal{F})$ and
$\mathcal{D}(M,\mathcal{F})$ that verifies $[{\mathcal
Q}(\mathcal{F})(\nabla(\mathcal{F}))(S)]=S,$ for all
$\nabla(\mathcal{F})\in\mathcal{C}(M,\mathcal{F})$ and all
$S\in\mathcal{S}(M,\mathcal{F})$, \item ${\mathcal
Q}(\mathcal{F})(\nabla(\mathcal{F}))={\mathcal
Q}(\mathcal{F})(\nabla'(\mathcal{F}))$, if $\nabla(\mathcal{F})$
and $\nabla'(\mathcal{F})$ are projectively \linebreak equivalent,
\item If $\Phi:(M,\mathcal{F})\to(M',\mathcal{F}')$ is a foliated
local diffeomorphism between two foliated manifolds
$(M,\mathcal{F})$ and $(M',\mathcal{F}')$, then
$${\mathcal Q}(\mathcal{F})(\Phi_{\mathcal C}^{*}\nabla(\mathcal{F'}))(\Phi_{\mathcal
S}^{*}S)(\Phi^{*}f)=\Phi^{*}({\mathcal
Q}(\mathcal{F}')(\nabla(\mathcal{F}'))(S)(f)),$$
 for all $\nabla(\mathcal{F}')\in\mathcal{C}(M',\mathcal{F}')$,
$S\in\mathcal{S}(M',\mathcal{F}')$, $f\in C^{\infty}(M',\mathcal{F}')$.
\end{itemize}
\end{defi}

Existence of a foliated natural and projectively invariant
quantization was proven in \cite{PoRaWo}.

\section{Singular objects}

Recall first the definition of a Riemannian orbifold.

\begin{defi}\label{RiemOrbi}
An $n$-dimensional ($n\in\mathbb{N}$; smooth or, more precisely,
$\Ci$-smooth) {\sl Riemannian orbifold structure} $V$ on a second
countable Hausdorff space $\vert V \vert$ is given by the
following data:
\begin{itemize}
\item An open cover $\lbrace V_{i}\rbrace_{i}$ of $\vert V\vert$.
\item For each $i\in I$, a connected and open subset
$U_{i}\subset\mathbb{R}^{n}$ with a Riemannian metric $h_i$; a
finite subgroup $\Gamma_{i}$ of isometries of the Riemannian
manifold $(U_{i},h_i)$; an open map $q_{i}:U_{i}\to V_{i}$, called
a local uniformization, that induces a homeomorphism from
$U_{i}/\Gamma_{i}$ onto $V_{i}$. \item For all $x_{i}\in U_{i}$
and $x_{j}\in U_{j}$ such that $q_{i}(x_{i})=q_{j}(x_{j})$, there
exist $W_{i}\subset U_{i}$ and $W_{j}\subset U_{j}$, open
connected neighborhoods of $x_{i}$ and $x_j$ respectively, and an
isometry $\phi_{ji}:W_{i}\to W_{j}$, called a change of charts,
such that $q_{j}\phi_{ji}=q_{i}$ on $W_{i}$.
\end{itemize}
\end{defi}

The assumption that the considered smooth orbifold be endowed with
a Riemannian metric is not a restriction, since any smooth
orbifold admits such a metric. Note further that any open subset
$U$ of any $n$-dimensional Riemannian orbifold, which is defined
by an orbifold atlas $\{(U_i,\Gamma_i,q_i)\}_i$, carries an induced
$n$-dimensional Riemannian orbifold structure defined by the atlas
$\{(\Omega_i:=q_i^{-1}(U\cap V_i),\Gamma_i,q_i|_{\Omega_i})\}_i$.

\begin{defi}
Let $f:V\to V'$ be a continuous map between two orbifolds $V$ and
$V'$. If for any $x\in V$, there exists a chart
$(U_{i},\Gamma_i,q_i)$ around $x$, i.e. such that $x\in
V_i=q_i(U_i)$, a chart $(U'_{j},\Gamma'_j,q'_j)$ around $f(x)$, as
well as a function $\tilde f\in C^{\infty}(U_{i},U'_{j})$, such
that $f q_{i}=q'_{j} \tilde f$, we say that $f$ is a {\sl smooth
map}. We denote by $C^{\infty}(V,V')$ the set of smooth mappings
from $V$ to $V'$ and by $\mathrm{Diff}(V,V')$ the set of
diffeomorphisms between $V$ and $V'$.
\end{defi}

In particular, a (continuous) function $f:V\to\R$ of an orbifold
$V$ is smooth, if for any $x\in V$, there is a chart
$(U_i,\Gamma_i,q_i)$ around $x$, such that $f q_i\in\Ci(U_i)$. If $U$
denotes an open subset of $V$, a (continuous) map $f:U\to\R$ is
smooth, if, for any $x\in U$, there exists a chart
$(U_i,\Gamma_i,q_i)$ in the neighborhood of $x$, such that $f
q_{i}\in C^{\infty}(q_{i}^{-1}(U\cap V_{i}))$. In the following
$C^{\infty}(U)$ denotes the associative commutative algebra of
smooth functions on $U$.

The assumption that $f:V\to \R$ be continuous is redundant here.
Indeed, since $q_i$ is surjective, we have $q_i(q_i^{-1}S_i)=S_i$,
for any $S_i\subset V_i$. Further, for any open $I\subset\R$, the
preimage $q_i^{-1}f|_{V_i}^{-1}I=(fq_i)^{-1}I$ is open and thus
$f|_{V_i}^{-1}I=q_i(q_i^{-1}f|_{V_i}^{-1}I)$ is open. Eventually,
we get $f^{-1}I=\cup_i f|_{V_i}^{-1}I$ is open in $V$.

\begin{defi} A {\sl differential operator} $D$ of order $k\ge 0$ of an orbifold $V$ is
an endomorphism of $C^{\infty}(V)$, such that we have on all
$U_i$,
$$(Df)q_{i}=\sum_{|\alpha|\le k}D_{\alpha}\partial_{x}^{\alpha}(f q_{i}),$$  where
$D_{\alpha}\in C^{\infty}(U_i)$ and where $k$ is independent of the
considered chart.\end{defi}

We denote by $\mathcal{D}^{k}(V)$ the $C^{\infty}(V)$-module of
differential operators of order $k$ of $V$ and by
$\mathcal{D}(V):=\bigcup_{i=0}^{\infty}\mathcal{D}^{i}(V)$ the Lie
algebra of all differential operators of $V$. As usual,
$[\mathcal{D}^i(V),\mathcal{D}^j(V)]\subset\mathcal{D}^{i+j-1}(V)$,
so that $\mathcal{D}^1(V)$ is a Lie subalgebra of $\mathcal{D}(V)$
and $\Ci(V)=\mathcal{D}^0(V)\subset\mathcal{D}^1(V)$ is a Lie
ideal of $\mathcal{D}^1(V)$.

\begin{defi}
The {\sl module of symbols} of degree $k\ge 0$ of $V$, which we
denote by $\mathcal{S}^{k}(V)$, is equal to
$\mathcal{D}^{k}(V)/\mathcal{D}^{k-1}(V)$. The module
$\mathcal{S}(V)$ of all symbols of $V$ is then equal
to $\bigoplus_{i=0}^{\infty}\mathcal{S}^{i}(V)$.
\end{defi}

\begin{defi}
The {\sl module and Lie algebra of vector fields} of $V$ is given by
$\mathrm{Vect}(V):=\mathcal{S}^{1}(V)$.
\end{defi}

\noindent{\bf Remark. }\rm{The map}
$\psi:\mathrm{Vect}(V)\ni[D]\mapsto D-D1\in\mathcal{D}^1(V)$ is a
splitting of the short exact sequence $$0\to \Ci(V)\to
\mathcal{D}^1(V)\to \mathrm{Vect}(V)\to 0$$ of $\Ci(V)$-modules,
so that
$$\mathcal{D}^1(V)\simeq \Ci(V)\oplus\mathrm{Vect}(V).$$ 

\begin{defi}
A {\sl torsion-free connection} $\nabla$ of $V$ is a bilinear map
$$\nabla : \mathrm{Vect}(V)\times\mathrm{Vect}(V)\to \mathrm{Vect}(V),$$
such that
\begin{itemize}
\item $\nabla_{fX}Y=f\nabla_{X}Y$, \item
$\nabla_{X}fY=(Xf)Y+f\nabla_{X}Y$, \item
$\nabla_{X}Y-\nabla_{Y}X=[X,Y]$,
\end{itemize}
for all $X\in\mathrm{Vect}(V)$, $Y\in\mathrm{Vect}(V)$ and $f\in C^{\infty}(V)$.
\end{defi}
We denote by $\mathcal{C}(V)$ the affine subspace of the space of
bilinear maps of $\mathrm{Vect}(V)$ that is made up by all
torsion-free connections of $V$.

\begin{defi}
A {\sl differential one-form} $\alpha$ of $V$ is a linear map from
$\mathrm{Vect}(V)$ to $C^{\infty}(V)$, such that for all
$X\in\mathrm{Vect}(V)$, we have on $V_i$,
$\alpha(X)\circ q_i=\sum_j\alpha_{j}X^{j}$, where
$X=[\sum_j X^{j}\partial_{x^j}]$ and $\alpha_{j}\in C^{\infty}(U_{i})$.
We denote by $\Omega^1(V)$ the $\Ci(V)$-module of
differential one-forms of $V$. 
\end{defi}

\begin{defi}
Two torsion-free connections $\nabla$ and $\nabla'$ of $V$ are
{\sl projectively equivalent} if and only if, for all vector
fields $X,Y\in\mathrm{Vect}(V),$
$$\nabla'_{X}Y=\nabla_{X}Y+\alpha(X)Y+\alpha(Y)X,$$
for some one-form $\alpha$ of $V$.
\end{defi}

\begin{defi}\label{local} A {\sl local isometry}
between two Riemannian orbifolds $V$ and $V'$ is a smooth map
$\varphi\in\Ci(V,V')$, such that for all $x\in V$, there exists a
chart $(U_i,\Gamma_i,q_i)$ of $V$, $x\in V_i:=q_i(U_i)$, and a chart
$(U'_j,\Gamma'_j,q'_j)$ of $V'$, $V'_j:=q_j'(U_j')$, such that
$\varphi\in\mathrm{Diff}(V_i,V'_j)$ admits a lift $\tilde\varphi:U_i\to
U'_j$, $\varphi q_i=q'_j\tilde\varphi$, which is an isometry between the
Riemannian manifolds $(U_i,h_i)$ and $(U'_j,h'_j)$, see Definition
\ref{RiemOrbi}.
\end{defi}

In the following definitions $\varphi$ denotes a local isometry
between two Riemannian orbifolds $V$ and $V'$ and notations are those of Definition \ref{local} (possible
extensions of these definitions are irrelevant for this paper).

\begin{defi}
The {\sl pullback of a function} $f\in C^{\infty}(V'_j)$ is
defined by $\varphi^*f:=f\circ\varphi\in C^{\infty}(V_i)$.
\end{defi}

\begin{defi}
The {\sl pullback of a $k$th order differential operator}
$D\in\mathcal{D}^k(V'_j)$,
$\varphi_{\mathcal{D}}^{*}D\in\mathcal{D}^k(V_i)$, is defined by
$$(\varphi_{\mathcal{D}}^{*}D)f:=\varphi^{*}(D(\varphi^{-1*}f)),$$
for all $f\in C^{\infty}(V_i)$.
\end{defi}

Indeed, we have
$\varphi_{\mathcal{D}}^{*}D\in\mathrm{End}(\Ci(V_i))$ and, since
$(U_i,\Gamma_i,\varphi q_i)$ is a compatible orbifold chart of $V'_j$,
we get on $U_i$
$$((\varphi_{\mathcal{D}}^{*}D)f)q_i=(D(f\varphi^{-1}))\varphi q_i=\sum_{|\alpha|\le
k}D_{\alpha}\partial_{x}^{\alpha}(fq_i),$$ with $D_{\alpha}\in\Ci(U_i)$.

It is easily checked that $\varphi^*_{\mathcal{D}}$ is a Lie algebra
isomorphism between $\mathcal{D}(V'_j)$ and
$\mathcal{D}(V_i)$.\medskip

Thanks to the fact that $\varphi_{\mathcal{D}}^{*}$ preserves the
order of the differential operators, one can give the following
definition:
\begin{defi}
If $S\in\mathcal{S}^k(V'_j)$ and if $S=[D]$ with
$D\in\mathcal{D}^k(V'_j)$, we define the {\sl symbol pullback} of
$S$ by
$$\varphi_{\mathcal{S}}^{*}S:=[\varphi_{\mathcal{D}}^{*}D]\in\mathcal{S}^k(V_i).$$
\end{defi}

\begin{defi}
The {\sl pullback map of vector fields} is
$$\varphi_{\mathrm{Vect}}^{*}:=\varphi_{\mathcal{S}}^{*}\vert_{\mathrm{Vect}(V'_j)}.$$
\end{defi}

Note first that if we identify the $\Ci(V)$-module
$\mathrm{Vect}(V)$ with the submodule
${\tt{Vect}}(V):=\psi(\mathrm{Vect}(V))$ of $\mathcal{D}^1(V)$,
see above, we have
\begin{equation}\varphi_{\mathrm{Vect}}^*=\varphi^*_{\mathcal
D}\vert_{{\tt{Vect}}(V_j')}.\label{PullVect2}\end{equation}

It follows immediately from the preceding definitions and the Lie
algebra isomorphism property of $\varphi^*_{\mathcal{D}}$ that
$\varphi^*_{\mathrm{Vect}}$ is a Lie algebra isomorphism between
$\mathrm{Vect}(V_j')$ and $\mathrm{Vect}(V_i)$. Further, for any
$f\in\Ci(V_j')$ and any $X\in\mathrm{Vect}(V'_j)$, we have
$$\varphi^*_{\mathrm{Vect}}(fX)=(\varphi^*f)(\varphi_{\mathrm{Vect}}^*X),$$
and, in view of Equation (\ref{PullVect2}), we also get
$$(\varphi^*_{\mathrm{Vect}}X)(\varphi^*f)=\varphi^*(Xf).$$

\begin{defi}
The {\sl pullback map of torsion-free connections}
$\varphi_{\mathcal{C}}^{*}:\mathcal{C}(V'_j)\to\mathcal{C}(V_i)$
is defined in this way:
$$(\varphi_{\mathcal{C}}^{*}\nabla)_{X}Y:=\varphi_{\rm{Vect}}^{*}(\nabla_{\varphi_{\rm{Vect}}^{-1*}X}\varphi_{\rm{Vect}}^{-1*}Y),$$
for all $\nabla\in\mathcal{C}(V_j')$, $X,Y\in\mathrm{Vect}(V_i)$.
\end{defi}

Remark that the just defined pullback of a torsion-free connection
is again a torsion-free connection, due to the preceding
properties of the pullback map for vector fields.

\begin{defi}
A {\sl natural and projectively invariant quantization $Q$ of
orbifolds} associates to any Riemannian orbifold $V$ a map
$$Q_{V}:\mathcal{C}(V)\times\mathcal{S}(V)\to\mathcal{D}(V),$$
such that
\begin{itemize}
\item $Q_{V}(\nabla)$ is a linear bijection between
$\mathcal{S}(V)$ and $\mathcal{D}(V)$, such that
$$[Q_{V}(\nabla)(S)]=S,$$ for all $\nabla\in\mathcal{C}(V)$ and all
$S\in\mathcal{S}^k(V)$, \item $Q_{V}(\nabla)=Q_{V}(\nabla')$, if
$\nabla$ and $\nabla'$ are projectively equivalent, \item if
$\varphi:V\to V'$ is a local isometry between two Riemannian
orbifolds $V$ and $V'$, then
$$Q_{V_i}(\varphi_{\mathcal{C}}^{*}\nabla)(\varphi_{\mathcal{S}}^{*}S)(\varphi^{*}f)=\varphi^{*}\left(Q_{V'_j}(\nabla)(S)(f)\right),$$
for all $\nabla\in\mathcal{C}(V'_j)$, $S\in\mathcal{S}(V'_j)$,
$f\in C^{\infty}(V'_j)$.
\end{itemize}
\end{defi}

\section{Resolution of a Riemannian orbifold}

For any $n$-dimensional Riemannian orbifold $V$, it is possible to
build a foliated manifold $\tilde{V}$, whose leaf space can be
identified with $V$. This construction is explained in details
e.g. in \cite{Girbau}. Let us briefly recall it here.\medskip

For any local uniformization $q_{i}:U_{i}\to V_{i}$, we denote by
$\tilde{U}_{i}(U_{i},\pi_i,O(n))$, where $O(n)$ is the orthogonal
group of degree $n$, the principal bundle of orthonormal frames of
the Riemannian manifold $(U_i,h_i)$. The $\Gamma_{i}$-action on
$U_{i}$ lifts in an obvious way to $\tilde{U}_{i}$: if $\tilde
u_{i}=(\tilde u_{i,1},\ldots,\tilde u_{i,n})\in\tilde{U}_{i}$ is
an orthonormal frame over $x_i\in U_i$ and if $g_{i}\in\Gamma_{i}$
is an isometry of $(U_i,h_i)$, then $g_{i}\tilde
u_{i}:=(g_{i*}\tilde u_{i,1},\ldots,g_{i*}\tilde u_{i,n})$ is an
orthonormal frame over $g_ix_i\in U_i$. This lifted action is
free, since an isometry is characterized by its derivative at one
point (more precisely, the map that associates to any $g_i\in
\Gamma_i$ an element $g_i\in\mathrm{Aut}(\tilde U_i)$ of the
automorphism group of the fiber bundle $\tilde U_i$ is a group
monomorphism). The quotient $\tilde{V}_{i}:=\tilde{U}_{i}/\Gamma_i$
is an ordinary smooth manifold. Indeed, as $\Gamma_{i}$ is a
finite group, its action on $\tilde{U}_{i}$ is also properly
discontinuous.

Similarly, any change of charts $\phi_{ji}:W_{i}\subset U_{i}\to
W_{j}\subset U_{j}$ lifts to a fiber bundle isomorphism
$$\tilde{\phi}_{ji}:\tilde{W}_{i}\subset\tilde{U}_{i}\to\tilde{W}_{j}\subset\tilde{U}_{j},\;\tilde
w_{i}\mapsto (\phi_{ji*}\tilde w_{i,1},\ldots,\phi_{ji*}\tilde
w_{i,n}).$$ Define now a projection
$$p_{i}:\tilde{V}_{i}\to V_{i}:[\tilde{u}_{i}]\mapsto
q_{i}\,\pi_{i}\,\tilde{u}_{i},$$ where $[.]$ denotes of course a
class of the quotient $\tilde V_i$. It is obviously well-defined.
Our goal is to glue the $\tilde V_i$ by means of gluing
diffeomorphisms
$$\tilde{f}_{ji}:p_{i}^{-1}(V_{ji})\subset\tilde{V}_{i}\to
p_{j}^{-1}(V_{ji})\subset\tilde{V}_{j},$$ where $V_{ji}=V_j\cap
V_i$, which verify the usual cocycle condition. Let $[\tilde
u_i]\in p_i^{-1}(V_{ji})$. Choose a representative $\tilde u_i$
(resp. $g_i\tilde u_i$), as well as a change of charts
$\phi_{ji}:W_i\subset U_i\to W_j\subset U_j$ such that
$\pi_i\,\tilde u_i\in W_i$ (resp. $\phi'_{ji}: g_iW_i\subset U_i\to
W_j'\subset U_j$), and set $$\tilde f_{ji}[\tilde
u_i]=[\tilde\phi_{ji}\,\tilde u_i]\in\tilde V_j\;
(\mbox{resp.}\;\tilde f_{ji}[\tilde
u_i]=[\tilde\phi'_{ji}\,g_i\,\tilde u_i]\in\tilde V_j).$$ Observe
that \begin{equation}p_j[\tilde\phi_{ji}\tilde
u_i]=q_j\,\pi_j\,\tilde\phi_{ji}\,\tilde
u_i=q_j\,\phi_{ji}\,\pi_i\,\tilde u_i=q_i\,\pi_i\,\tilde
u_i=p_i[\tilde u_i]\in V_{ji},\label{GluingProj}\end{equation} and
that the map $\tilde f_{ji}$ is well-defined, since the two chart
changes $\phi_{ji}$ and $\phi_{ji}'\,g_i$ defined on $W_i$ coincide
up to $g_j\in\Gamma_j$. Eventually, it is well-known that the chart
changes $\phi_{ji}$ verify the cocycle equation
$g_{ijk}\phi_{ki}=\phi_{kj}\phi_{ji}$, $g_{ijk}\in\Gamma_k$; this
entails that the same equation holds true for the lifts $\tilde
\phi_{ji}$ and thus that we have $\tilde f_{ki}=\tilde f_{kj}\tilde
f_{ji}$. Hence, if we glue the $\tilde{V}_{i}$ according to the
$\tilde{f}_{ji}$, we get a smooth manifold $\tilde{V}$ of
dimension $n(n+1)/2$.

Let now $\tilde V_i\ni[\tilde u_i]\simeq [\tilde \phi_{ji}\tilde
u_i]\in \tilde V_j$ be an element of $\tilde V$. It follows from
Equation (\ref{GluingProj}) that the local projections $p_i:\tilde
V_i\to V_i$ define a global projection $p:\tilde{V}\to V$.
Moreover, the manifold $\tilde{V}$ admits a right $O(n)$-action.
Indeed, for any $i$, the canonical ``matrix product'' right action
of $M\in O(n)$ on an orthonormal frame $\tilde u_i\in \tilde U_i$
is an orthonormal frame over the same point. Since clearly
$(g_i\tilde u_i)M=g_i(\tilde u_iM)$, this $O(n)$-action on $\tilde
U_i$ induces an action on $\tilde V_i$, given by
$[\tilde{u}_{i}]M:=[\tilde{u}_{i}M]$. Thanks to the fact that we
also have $(\tilde{\phi}_{ji}\tilde u_i)M=\tilde\phi_{ji}(\tilde
u_iM)$, we get a global $O(n)$-action on $\tilde V$. The orbits of
this action, which coincide with the fibers of the projection
$p:\tilde{V}\to V$, are known to be the leaves of a regular
foliation $\mathcal{F}$ on $\tilde{V}$.

We can find an atlas of $\tilde V$ made up by charts that are
adapted to $\mathcal F$. It suffices to build such an atlas for
$\tilde V_i=\tilde U_i/\Gamma_i$ by means of the general technique
for quotients of manifolds by free and properly discontinuous
group actions. Let $[\tilde{u}_{i}]\in\tilde{V}_{i}$ and let
$\tilde{U}$ be a neighborhood of $\tilde{u}_{i}$ in
$\tilde{U}_{i}$ such that $g_{i}\tilde{U}\cap\tilde{U}=\emptyset$,
for all $g_{i}\in\Gamma_{i}$ different from the identity. Such a
neighborhood exists since the action of $\Gamma_{i}$ is properly
discontinuous. We may assume that $\tilde U$ is contained in an open of
trivialization. For any $[\tilde{u}]\in[\tilde{U}]$, there is a unique
representative, say $\tilde u$, in $\tilde U.$ The
coordinates of $[\tilde u]$ are then $(M_{\tilde u},\pi_{i}\,\tilde{u})$, where
$M_{\tilde u}\in O(n)$ is the orthogonal matrix associated to $\tilde u$ via the
trivialization. It is a matter of
common knowledge that the coordinate systems
$$\psi:[\tilde U]\ni[\tilde u]\mapsto (M_{\tilde u},\pi_i\,\tilde
u)\in O(n)\times\pi_{i}\tilde{U}$$ form an atlas of $\tilde V_i$. Further, they
are obviously adapted to $\mathcal{F}$, the transverse coordinates
of $[\tilde{u}]$ being the components of
$\pi_{i}\,\tilde{u}$.\medskip

Observe that $p[\tilde U]=q_i\pi_i\tilde U$ is an open subset
of the orbifold $V_i$ defined by the chart $(U_i,\Gamma_i,q_i)$, so
that it is itself an orbifold for the chart
$(\Omega_i:=q_i^{-1}(q_i\pi_{i}\tilde{U}),\Gamma_i,$ $q_i|_{\Omega_i})$. 

\section{Singular quantization}

In the following $V$ denotes a Riemannian orbifold and $(\tilde
V,\mathcal{F})$ is its foliated resolution.

\begin{prop}
The map
$$p^*:C^{\infty}(V)\ni f\mapsto f\,p\in C^{\infty}(\tilde{V},\mathcal{F})$$
is a linear isomorphism.
\end{prop}

\begin{proof} If $f\in\Ci(V)$, it is clear that $f\,p$ is a foliated function.
Since around any point of $\tilde{V}$ there is a chart $([\tilde
U],\psi)$, such that
$$f\,p\,\psi^{-1}=f\,q_i\,\pi_i\,\psi^{-1}=f\,q_i\,\mathrm{pr}_2,$$
where $\mathrm{pr}_2$ is the projection from $O(n)\times U$ onto
$U$, the function $f\,p$ is also smooth. Conversely, a foliated
function gives rise to a function of the leaf space, i.e. to a
function of $V$.\end{proof}

\begin{prop}
The map
$$p_{\mathcal{D}}^{*}:\mathcal{D}^k(V)\ni D\mapsto p^{*}\,D\,p^{*-1}\in
\mathcal{D}^k(\tilde{V},\mathcal{F})$$
is a linear isomorphism and even a Lie algebra isomorphism between
$\mathcal{D}(V)$ and $\mathcal{D}(\tilde{V},\mathcal{F})$.
\end{prop}

\begin{proof} The unique point that requires an explanation is the fact that the
conjugate operator has the appropriate local form. This question
is actually just a matter of notations. Observe that if the
variable $[\tilde u]$ runs through an adapted chart domain
$[\tilde U]\subset \tilde V_i$, then $p[\tilde u]=q_i\pi_i\tilde
u=:q_i y,$ where the transverse coordinates $y=(y^1,\ldots,y^n)$
run through the corresponding open subset $U\subset U_i$. Further,
as aforementioned, a foliated function $g[\tilde u]$ factors in
the form $\tilde g p[\tilde u]=\tilde g\, q_i y,$ where $\tilde g$
is a singular function. Hence, if $D\in\mathcal{D}^k(V)$, the
value at $g$ of the endomorphism $p^*_{\mathcal{D}}D =
p^*\,D\,p^{*-1}$ locally reads

$$\begin{array}{c}(p^*_{\mathcal{D}}D)(g)[\tilde u]=D(\tilde g)\,p[\tilde u]=D(\tilde
g)\,q_iy=\sum_{|\alpha|\le k}\tilde D_{\alpha}y\;
\partial_y^{\alpha}(\tilde g\,q_iy) \\\\=\sum_{|\alpha|\le
k}\tilde D_{\alpha}\pi_i\tilde u\,
\partial_y^{\alpha}(\tilde gp[\tilde u])=\sum_{|\alpha|\le
k}\tilde D_{\alpha}pr_{2}(M,y)\;
\partial_y^{\alpha}(g(M,y)),\end{array}$$

\vspace{5mm}\noindent where we identified the point $[\tilde u]$
with its coordinates $(M,y)$.\end{proof}

\begin{prop}
The map
$$p_{\mathcal{S}}^{*}:\mathcal{S}^k(V)\ni
[D]\mapsto[p_{\mathcal{D}}^{*}D]\in\mathcal{S}^k(\tilde{V},\mathcal{F})$$
is a linear isomorphism.
\end{prop}

\begin{proof}
Obvious.
\end{proof}

The restriction of the mapping $p_{\mathcal{S}}^{*}$ to
$\mathcal{S}^1(V)$ is of course a Lie algebra isomorphism
$p^*_{\mathrm{Vect}}$ between $\mathrm{Vect}(V)$ and
$\mathrm{Vect}(\tilde{V},\mathcal{F})$. Furthermore, just as for
the pullback by a local isometry, we have
$p^*_{\mathrm{Vect}}(fX)=(p^*f)(p^*_{\mathrm{Vect}}X)$ and
$(p^*_{\mathrm{Vect}}X)(p^*f)=p^*(Xf)$, for all $f\in\Ci(V)$ and
all $X\in\mathrm{Vect}(V)$.\medskip

\noindent{\bf Remark :} One can easily show that the previous results can be
extended to the case where $\tilde{V}$ is replaced by an open set $\tilde{\Omega}$
of $\tilde{V}$ and where $V$ is replaced by $p(\tilde{\Omega})$.

\begin{lem}
There exists a pullback $p_{\Omega}^{*}$ that maps singular
1-forms of $V$ to foliated 1-forms of $(\tilde{V},\mathcal F)$ and
verifies
$$(p_{\Omega}^{*}\alpha)(X)=p^{*}(\alpha(p_{\mathrm{Vect}}^{*-1}(X)),$$
for all $\alpha\in\Omega^1(V)$ and all
$X\in\mathrm{Vect}(\tilde{V},\mathcal{F})$.
\end{lem}

\begin{proof}
Let $\alpha\in\Omega^1(V)$ and $X\in\mathrm{Vect}(\tilde{V})$.
Note that for the moment we do not assume that $X$ is foliated.
For any chart $([\tilde{U}],(M,y))$ of $\tilde V$ adapted to
$\mathcal{F}$, we can apply the preceding pullback results to the
orbifold $p[\tilde U]$. If $X$ reads $X=\sum_{\iota}
X^{\iota}\partial_{M^{\iota}}+\sum_{\frak i} X^{\frak
i}\partial_{y^{\frak i}}$ in $[\tilde U]$, we thus can set
$$(p_{\Omega}^{*}\alpha)(X)\vert_{[\tilde{U}]}:=\sum_{\frak i}X^{\frak
i}p^{*}(\alpha(p_{\mathrm{Vect}}^{*-1}[\partial_{y^{\frak i}}]))\in\Ci([\tilde
U]),$$
where the second factors of the {\small RHS} are foliated locally
defined functions. One can quite easily prove that the functions
$(p_{\Omega}^{*}\alpha)(X)\vert_{[\tilde{U}]}$ can be glued and
yield a global function $(p_{\Omega}^{*}\alpha)(X)$ of $\tilde V$,
since, if $(N,z)$ are other adapted coordinates, we have $z=z(y)$.
It follows that $p^*_{\Omega}\alpha$ is a differential 1-form of $\tilde
V$, which is clearly foliated in view of the preceding definition.
Observe eventually that for foliated vector fields $X$, the
{\small RHS} of the defining equation reads
$$p^{*}(\alpha(p_{\mathrm{Vect}}^{*-1}\,[\sum_{\frak i}X^{\frak
i}\partial_{y^{\frak i}}])).$$\end{proof}


\begin{prop}
The map
$$p_{\mathcal{C}}^{*}:\mathcal{C}(V)\ni\nabla\mapsto
p_{\mathcal{C}}^{*}\nabla\in\mathcal{C}(\tilde{V},\mathcal{F}),$$
where $p_{\mathcal{C}}^{*}\nabla$ is defined by
$$(p_{\mathcal{C}}^{*}\nabla)_{X}Y=p_{\mathrm{Vect}}^{*}(\nabla_{p_{\mathrm{Vect}}^{*-1}X}p_{\mathrm{Vect}}^{*-1}Y),$$
transforms projective classes of singular torsion-free connections
in projective classes of foliated torsion-free connections.
\end{prop}

\begin{proof}
The result is a consequence of the preceding propositions.
\end{proof}
\begin{thm}
There exists a natural and projectively invariant quantization of
orbifolds. If $Q$ denotes this quantization and if $V$ is a
Riemannian orbifold, the map $Q_{V}$ is, for any singular
connection $\nabla\in\mathcal{C}(V)$ and any singular symbol
$S\in\mathcal{S}(V)$, defined by:
$$Q_{V}(\nabla)(S):=p_{\mathcal D}^{*-1}\left({\mathcal
Q}(\mathcal{F})(p_{\mathcal{C}}^{*}\nabla)(p_{\mathcal{S}}^{*}S)\right),$$
where ${\mathcal Q}(\mathcal F)$ is the map associated by the
foliated natural and projectively invariant quantization $\mathcal
Q$ to the foliated manifold $(\tilde V,\mathcal F)$.
\end{thm}
\begin{proof} The unique required property of $Q$, which is not
obvious in view of the above propositions and of the properties of
$\mathcal Q$, is its naturality.

Let $\varphi:V\to V'$ be a local isometry between two Riemannian
orbifolds $V,V'$ and let $\tilde\varphi:U_i\to U_j'$ be the isometry
that lifts the diffeomorphism $\varphi:V_i\to V_j'$. Then
$\tilde\varphi_*:\tilde U_i\to \tilde U_j'$ is a bundle isomorphism
over $\tilde \varphi$, which, in view of standard arguments, induces
a diffeomorphism $\Phi:\tilde V_i\to\tilde V_j'$, $\Phi[\tilde
u_i]=[\tilde\varphi_*\tilde u_i]$. It follows that
$$p'\Phi[\tilde u_i]=q'_j\pi'_j\tilde\varphi_*\tilde
u_i=q'_j\tilde\varphi\pi_i\tilde u_i=\varphi q_i\pi_i\tilde u_i=\varphi
p[\tilde u_i],$$ so that \begin{equation}\label{Comm} p'\Phi=\varphi
p,\end{equation} where notations are self-explaining. Further,
$\Phi:\tilde V_i\to \tilde V_j'$ is a foliated local diffeomorphism
between $(\tilde V_i,\mathcal F)$ and $(\tilde V_j',\mathcal F')$.
Indeed, it maps any leaf $p^{-1}v_i$, $v_i\in V_i$ of $\mathcal F$
into a leaf of $\mathcal F'$, since $p'\Phi p^{-1}v_i=\varphi
pp^{-1}v_i=\{\varphi v_i\}.$

It is straightforwardly checked that equation (\ref{Comm}) entails
\begin{equation}p^*\varphi^*=\Phi^*p'^*,\,\; p^*_{\mathcal
S}\,\varphi_{\mathcal S}^*=\Phi^*_{\mathcal S}\,p'^*_{\mathcal S},\;\,
p^*_{\mathcal C}\,\varphi^*_{\mathcal C}=\Phi^*_{\mathcal
C}\,p'^*_{\mathcal C}.\label{Comm*}\end{equation} The definition
of the singular quantization (which implies a similar equation for
$Q_V(\nabla)(S)(f)$), the commutation relations (\ref{Comm*}), and
the naturality of the foliated quantization finally show that the
singular quantization is natural as well.

\end{proof}
\section{Acknowledgements}
The research of Norbert Poncin was supported by UL-grant
SGQnp2008. Fabian Radoux thanks the Belgian FNRS for his Research
Fellowship.\\

\end{document}